\input amstex
\documentstyle{amsppt}
\document
\topmatter
\title
Semi-symmetric K\"ahler surfaces.
\endtitle
\author
W\l odzimierz Jelonek
\endauthor

\abstract{The aim of this paper is to describe   K\"ahler surfaces
which admit an opposite almost Hermitian structure satisfying the
first Gray condition. }
\thanks{MS Classification: 53C55,53C25,53B35. Key words and phrases:
K\"ahler surface, holomorphic sectional curvature, quasi constant
holomorphic sectional curvature,QCH manifold, ambik\"ahler
mani\-fold }\endthanks
 \endabstract
\endtopmatter
\define\G{\Gamma}
\define\bJ{\overline{J}}
\define\DE{\Cal D^{\perp}}
\define\e{\epsilon}
\define\n{\nabla}
\define\om{\omega}
\define\Om{\Omega}
\define\r{\rightarrow}
\define\w{\wedge}
\define\k{\diamondsuit}
\define\th{\theta}

\define\a{\alpha}

\define\lb{\lambda}

\define\1{D_{\lb}}
\define\2{D_{\mu}}
\define\0{\Omega}

\define\De{\Cal D}

\define\m{(M,g,J)}
\define \E{\Cal E}
\bigskip
{\bf 0. Introduction. } The aim of the present paper is to
describe connected K\"ahler surfaces $(M,g,J)$ admitting a
negative almost hermitian structure $\bJ$ satisfying the first
Gray condition  $R(\bJ X,\bJ Y, Z, U)=R(X,Y,Z,U)$.

Such surfaces are QCH K\"ahler surfaces (see [J-4]) i.e. surfaces
admitting a global, $2$-dimensional, $J$-invariant distribution
$\De$ having the following property: The holomorphic curvature
$K(\pi)=R(X,JX,JX,X)$ of any $J$-invariant $2$-plane $\pi\subset
T_xM$, where $X\in \pi$ and $g(X,X)=1$, depends only on the point
$x$ and the number $|X_{\De}|=\sqrt{g(X_{\De},X_{\De})}$, where
$X_{\De}$ is an orthogonal projection of $X$ on $\De$. In this
case  we have
$$R(X,JX,JX,X)=\phi(x,|X_{\De}|)$$ where $\phi(x,t)=a(x)+b(x)t^2+c(x)t^4$ and
 $a,b,c$ are smooth functions on $M$. Also $R=a\Pi+b\Phi+c\Psi$
 for certain curvature tensors $\Pi,\Phi,\Psi\in \bigotimes^4\frak X^*(M)$
  of K\"ahler type. The investigation of such manifolds, called QCH K\"ahler manifolds,  was
started by G. Ganchev and V. Mihova in [G-M-1],[G-M-2]. Every QCH
K\"ahler surface is holomorphically pseudosymmetric and
$R.R=\frac16(\tau-\kappa)\Pi.R$ (see [J-4],[O]). In our paper
[J-2] we used their local results to obtain a global
classification of such manifolds under the assumption that $\dim
M=2n\ge 6$.
  In
 the present paper we show that a K\"ahler surface $\m$ is
 semi-symmetric
 if and only if is locally symmetric or it admits a negative almost
 hermitian structure $\bJ$ which satisfies the first Gray
 condition $R(\bJ X,\bJ Y, Z, U)=R(X,Y,Z,U)$.
 We also prove that  a semmi-symetric K\"ahler surface $\m$ is a QCH K\"ahler surface
 or $\m$ is locally isometric to a space form.  In [J-4] we have proved that
 $(M,g,J)$ is a QCH K\"ahler surface iff it admits a negative
 almost complex structure $ \overline{J}$ satisfying the Gray
 second condition
 $R(X,Y,Z,W)-R(\overline{J}X,\overline{J}Y,Z,W)=R(\overline{J}X,Y,\overline{J}Z,W)+R(\overline{J}X,Y,Z,\overline{J}W)$.
  In [A-C-G]  Apostolov, Calderbank and Gauduchon have classified weakly selfdual K\"ahler surfaces, extending
  the result of Bryant who classified self-dual K\"ahler surfaces [B].
  Weakly self-dual K\"ahler surfaces turned out to be of Calabi
  type and of orthotoric type or surfaces with parallel Ricci tensor.
   Any Calabi type K\"ahler surface and every
 orthotoric K\"ahler surface is a QCH manifold. In both cases the
 opposite complex strucure $\bJ$ is conformally K\"ahler.
\bigskip
{\bf 1. The first Gray condition.} Let $\m$ be a $4$-dimensional
K\"ahler manifold with a negative almost Hermitian structure
$\bJ$. Then $\De=ker(J\bJ-Id)$ is a $J$-invariant distribution .
Let $\frak X(M)$ denote the algebra of all differentiable vector
fields on $M$. If $X\in\frak X(M)$ then by $X^{\flat}$ we shall
denote the 1-form $\phi\in \frak X^*(M)$ dual to $X$ with respect
to $g$, i.e. $\phi(Y)=X^{\flat}(Y)=g(X,Y)$. By $\om$ we shall
denote the K\"ahler form of $\m$ i.e. $\om(X,Y)=g(JX,Y)$. Let $\m$
be a QCH K\"ahler surface with respect to $J-invariant$
2-dimensional distribution $\De$. Let us denote by $\E$ the
distribution $\DE$, which is a $2$-dimensional, $J$-invariant
distribution. Then $\E=ker(J\bJ+Id)$. By $h,m$ respectively we
shall denote the tensors $h=g\circ (p_{\De}\times
p_{\De}),m=g\circ (p_{\E}\times p_{\E})$, where $p_{\De},p_{\E}$
are the orthogonal projections on $\De,\E$ respectively. It
follows that $g=h+m$. For every almost Hermitian manifold $\m$ the
self-dual Weyl tensor $W^+$ decomposes under the action of the
unitary group $U(2)$. We have $\bigwedge^*M=\Bbb R\oplus LM$ where
$LM=[[\bigwedge^{(0,2)}M]]$ and we can write $W^+$ as a matrix
with respect to this block decomposition
$$W^+=\pmatrix \frac{\kappa}6& W_2^+\cr (W_2^+)^*&W_3^+-\frac{\kappa}{12}Id_{|LM}\endpmatrix$$
where $\kappa$ is the conformal scalar curvature of $(M,g,J)$. By
$\tau$ we denote a scalar curvature of $\m$. The selfdual Weyl
tensor $W^+$ of $\m$ is called degenerate if $W_2=0,W_3=0$. In
general the self-dual Weyl tensor of 4-manifold $(M,g)$ is called
degenerate if it has at most two eigenvalues as an endomorphism
$W^+:\bigwedge^+M\r\bigwedge^+M$. We say that an almost Hermitian
structure $J$ satisfies the first Gray condition if
$$R(X,Y,Z,W)=R(JX,JY,Z,W)\tag{G1}$$
and the second Gray curvature condition if
$$R(X,Y,Z,W)-R(JX,JY,Z,W)=R(JX,Y,JZ,W)+R(JX,Y,Z,JW),\tag {G2}$$
which is equivalent to $Ric(J,J)=Ric$ and $W_2^+=W_3^+=0$.
Condition $(G_1)$ implies $(G_2)$.  $\m$ satisfies the second Gray
condition if $J$ preserves the Ricci tensor and $W^+$ is
degenerate.  A K\"ahler surface is QCH if and only if it admits a
negative almost hermitian struce satisfying the second Gray
condition. Every QCH K\"ahler surface is holomorphically
pseudosymmetric and $R.R=\frac16(\tau-\kappa)\Pi.R$ (see
[J-4],[O]). We shall denote by $Ric_0$ and $\rho_0$ the trace free
part of the Ricci tensor $Ric$ and the Ricci form $\rho$
respectively. An ambik\"ahler structure on a real 4-manifold
consists of a pair of K\"ahler metrics $(g_+,J_+,\om_+)$ and
$(g_-,J_-,\om_-)$ such that $g_+$ and $g_-$ are conformal metrics
and $J_+$ gives an opposite orientation to that given by $J_-$
(i.e the volume elements $\frac12\om_+\w\om_+$ and
$\frac12\om_-\w\om_-$ have opposite signs).  A foliation $\Cal F$
on a Riemannian manifold $(M,g)$ is called conformal if
$$L_Vg=\a(V)g$$ holds on $T\Cal F^{\perp}$ where $\a$ is a one
form vanishing on $T\Cal F^{\perp}$.  A foliation $\Cal F$ is
called homothetic if is conformal and $d\a=0$  (see  [Ch-N] ). A
foliation $\Cal F$ on a complex manifold $(M,J)$ is called complex
if $J\Cal F\subset \Cal F$ and is called holomorphic if
$L_XJ(TM)\subset\Cal F$ for any $X\in\G(\Cal F)$. Complex
homothetic foliations by curves on K\"ahler manifolds were
  recently  classified locally in [Ch-N].
\medskip
{\bf 2. Curvature tensor of a QCH K\"ahler surface.} We shall
recall some results from [G-M-1]. Let
$$R(X,Y)Z=([\n_X,\n_Y]-\n_{[X,Y]})Z\tag 2.1$$ and let us write $$R(X
,Y,Z,W)=g(R(X,Y)Z,W).$$ If $R$ is the curvature tensor of a QCH
K\"ahler manifold $\m$, then there exist functions $a,b,c\in
C^{\infty}(M)$ such that
$$R=a\Pi+b\Phi+c\Psi,\tag 2.2$$
where $\Pi$ is the standard K\"ahler tensor of constant
holomorphic curvature i.e.
$$\gather \Pi(X,Y,Z,U)=\frac14(g(Y,Z)g(X,U)-g(X,Z)g(Y,U)\tag 2.3\\
+g(JY,Z)g(JX,U)-g(JX,Z)g(JY,U)-2g(JX,Y)g(JZ,U)),\endgather $$
the tensor $\Phi$ is defined by the following relation
$$\gather \Phi(X,Y,Z,U)=\frac18(g(Y,Z)h(X,U)-g(X,Z)h(Y,U)\tag 2.4\\+g(X,U)h(Y,Z)-g(Y,U)h(X,Z)
+g(JY,Z)h(JX,U)\\-g(JX,Z)h(JY,U)+g(JX,U)h(JY,Z)-g(JY,U)h(JX,Z)\\
-2g(JX,Y)h(JZ,U)-2g(JZ,U)h(JX,Y)),\endgather$$ and finally
$$\Psi(X,Y,Z,U)=-h(JX,Y)h(JZ,U)=-(h_J\otimes h_J)(X,Y,Z,U).\tag 2.5$$
where $h_J(X,Y)=h(JX,Y)$.
 Let $V=
 (V,g,J)$ be a real $2n$ dimensional vector space with
complex structure $J$ which is skew-symmetric with respect to the
scalar product $g$ on $V$.  Let assume further that  $V= D\oplus
E$ where $D$ is a 2-dimensional, $J$-invariant subspace of $V$,
$E$ denotes its orthogonal complement in $V$. Note that  the
tensors $\Pi,\Phi,\Psi$ given above are of K\"ahler type.  It is
easy to check that for a unit vector $X\in V$
$\Pi(X,JX,JX,X)=1,\Phi(X,JX,JX,X)=|X_{D}|^2,\Psi(X,JX,JX,X)=|X_{D}|^4$,
where $X_D$ means an orthogonal projection of a vector $X$ on the
subspace $D$ and $|X|=\sqrt{g(X,X)}$. It follows that for a tensor
$(2.2)$ defined on $V$ we have
$$R(X,JX,JX,X)=\phi(|X_D|)$$ where $\phi(t)=a+bt^2+ct^4$.

If $Ric_0=\delta(h-m)$ then (see [J-4])
$$R=(\frac{\tau}6-\delta+\frac\kappa{12})\Pi+(2\delta-\frac{\kappa}2)\Phi+\frac{\kappa}2\Psi.\tag 2.6$$
Let $J, \overline{J}$ be hermitian, opposite orthogonal structures
on a Riemannian 4-mani\-fold $(M,g)$ such that $J$ is a positive
almost complex structure. Let $\Cal
E=ker(J\overline{J}-Id),\De=ker(J\overline{J}+Id)$ and let the
tensors $\Pi,\Phi,\Psi$ be defined as above where
$h=g(p_{\De},p_{\De})$. Let us define a tensor
$K=\frac16\Pi-\Phi+\Psi$.  Then $K$ is a curvature tensor,
$b(K)=0,c(K)=0$ where $b$ is Bianchi operator and $c$ is the Ricci
contraction. For a QCH K\"ahler surface $\m$ we have
$W^-=\frac{\kappa}2K$ (see [J-4]).

\medskip
Let $\m$ be a K\"ahler surface which is a QCH manifold with
respect to the distribution $\De$.  Then $\m$ is also QCH manifold
with respect to the distribution $\Cal E=\De^{\perp}$ and if
$\Phi',\Psi'$ are the above tensors with respect to $\Cal E$ then
$$R=(a+b+c)\Pi-(b+2c)\Phi'+c\Psi'.\tag 2.7$$

\medskip
 If $\m$ is a QCH K\"ahler surface then one can show that
the Ricci tensor $\rho$ of $\m$ satisfies the equation
$$\rho(X,Y)=\lb m(X,Y)+\mu h(X,Y)\tag 2.8$$
where $\lb=\frac{3}2a+\frac  b4,\mu=\frac{3}2a+\frac{5}4b+c$ are
eigenvalues of $\rho$ (see [G-M-1], Corollary 2.1 and Remark 2.1.)
In particular the distributions $\E,\De$ are eigendistributions of
the tensor $\rho$ corresponding to the eigenvalues $\lb,\mu$ of
$\rho$.

{\bf Lemma 1.}  {\it  The tensors $\Pi.\Phi,\Pi.\Psi$ are linearly
independent.}

\medskip
{\it Proof.}  Let  $\{e_1,\e_2,\e_3,\e_4\}$ be an orthogonal basis
in $TM$ such that  $$\De=span\{e_1,e_2\},\E=span\{e_3,e_4\}$$ and
$e_2=Je_1,e_4=Je_3$.   Then

$$\gather
\Pi(e_1,e_3).\Phi(e_1,e_4,e_3,e_4)=\Phi(e_3,e_4,e_3,e_4)-\Phi(e_1,e_2,e_3,e_4)\tag 2.9\\
-\Phi(e_1,e_4,e_1,e_4)-\Phi(e_1,e_4,e_3,e_2)=\frac3{16}\\
\Pi(e_1,e_3).\Psi(e_1,e_4,e_3,e_4)=0\endgather$$ $\k$
\medskip
{\bf Proposition 1.} {\it If a QCH K\"ahler surface satisfies at a
point $x\in M$ the condition $R.R=0$ then at $x$ we have $R=a\Pi$
or $2a+b=\frac16(\tau-\kappa)=0$.}
\medskip
{\it Proof.}  We have $R.R=(a+\frac b2)\Pi.R$. On the other hand
$\Pi.R=b\Pi.\Phi+c\Pi.\Psi$. From the Lemma $\Pi.R=0$ if and only
if $b=c=0$ and consequently if $R=a\Pi$.  If $R\ne a\Pi$ then
$R.R=0$ implies $2a+b=0$.$\k$

\medskip
{\bf Lemma 2.}  {\it  The tensor $R=a\Pi+b\Phi+c\Psi$ satisfies
the first Gray condition with respect to $\bJ$ if and only if
$2a+b=0$.}
\medskip
{\it Proof} Note that $\Psi(\bJ X,\bJ Y, Z, U)=\Psi(X,Y,Z,U)$. We
also have ($I=\bJ$)
$$\gather \Pi(\bJ X,\bJ Y, Z,U)-\Pi(X,Y,Z,U)=\tag 2.10\\ \frac12(m(IY,Z)h(IX,U)+m(IX,U)h(IY,Z)-m(IX,Z)h(IY,U)\\
-h(IX,Z)m(IY,U)-h(Y,Z)m(X,V)-m(Y,Z)h(X,U)\\+h(X,Z)m(Y,U)+m(X,Z)h(Y,U))\endgather$$
and
$$\gather \Phi(\bJ X,\bJ Y, Z,U)-\Phi(X,Y,Z,U)=\tag 2.11\\ \frac14(m(IY,Z)h(IX,U)+m(IX,U)h(IY,Z)-m(IX,Z)h(IY,U)\\
-h(IX,Z)m(IY,U)-h(Y,Z)m(X,V)-m(Y,Z)h(X,U)\\+h(X,Z)m(Y,U)+m(X,Z)h(Y,U))\endgather$$

Hence

$$\gather R(\bJ X,\bJ Y, Z,U)-R(X,Y,Z,U)=\tag 2.12\\(2a+b)(m(IY,Z)h(IX,U)+m(IX,U)h(IY,Z)-m(IX,Z)h(IY,U)\\
-h(IX,Z)m(IY,U)-h(Y,Z)m(X,V)-m(Y,Z)h(X,U)+\\h(X,Z)m(Y,U)+m(X,Z)h(Y,U))\endgather$$

Consequently $R(\bJ X,\bJ Y, Z,U)-R(X,Y,Z,U)=0$ if and only if
$2a+b=0$.
\medskip
{\bf Proposition 2.} {\it  Let us assume that a K\"ahler $\m$
surface admits a negative almost Hermitian structure $\bJ$
satisfying the first Gray condition.  Then $\m$ is a QCH
semi-symmetric surface and $\tau=\kappa$ where $\kappa$ is the
conformal scalar curvature of $\bJ$.}
\medskip
{\it Proof.}   Since $\m$ satisfies the first Gray condition with
respect to $\bJ$ it clearly satisfies the second condition and is
a QCH surface. On the other hand $R.R=(a+\frac b2)\Pi.R=0$ since
$a+\frac b2=\frac16(\tau-\kappa)=0$.$\k$
\medskip
{\bf Lemma  3. }{\it Let us assume that a product $M=\Bbb R\times
N$ is a K\"ahler surface where $N$ is 3-dimensional Riemannian
manifold. Then $M$ is locally a product of Riemannian surfaces. If
$M$ is simply connected and complete then $M$ is a product of
Riemannian surfaces }
\medskip
{\it Proof.}  Let $H$ be a unit vector field tangent to $\Bbb R$.
Then $\n H=0$.  Thus if $X=JH$ then $X$ is a unit covariantly
constant vector field in $N$.  The distribution $\De=\{Y\in
TN:g(X,Z)=0\}$ is parallel and $J\De=\De$.  Hence locally $N=\Bbb
R\times \Sigma$ and $M$ is a product of Riemannian surfaces. If
$M$ is simply connected and complete then $N$ is simply connected
and complete and from the De Rham theorem $N=\Bbb R\times \Sigma$
and $M=\Bbb C\times \Sigma$. $\k$

\medskip
{\bf Proposition 4.} {\it Let $\m$ be a semi-symmetric K\"ahler
surface. Then locally $\m$ is a  space form or a QCH K\"ahler
surface.}

\medskip
{\it Proof.}  We use a classification result of Szabo [Sz] and
Lumiste [L] and Lemma 3. Note that $JV^0=V^0$ in the Szabo
decomposition since $R(X,Y)\circ J=J\circ R(X,Y)$. Hence $dim
V^0=0,2,4$.  Note also that for elliptic,hyperbolic and Euclidean
cone we have $dim V^0=1$. Hence locally $\m$ is a symmetric space,
a product of two Riemannian surfaces and a space foliated by two
dimensional Euclidean space. A space foliated by two dimensional
Euclidean space is a QCH K\"ahler surface with respect to $\E=V^0$
or $\De=V^1$. In fact
$R(X,JX,JX,X)=R(X_{\De},JX_{\De},JX_{\De},X_{\De})$ where
$X_{\De}$ is an orthogonal projection of $X$ onto $\De$. Hence
$R=c\Psi$ where $\Psi$ is the tensor with respect to $\De$. It is
also clear that product of Riemannian surfaces
$M=\Sigma_1\times\Sigma_2$ is a QCH K\"ahler surface with respect
to $\De=T\Sigma_1$ or $\E=T\Sigma_2$.   Note also that locally
symmetric irreducible K\"ahler surface is self-dual hence is a
space form  (see [D]).$\k$
\medskip
{\bf Proposition 5.  }{\it  Let us assume that $\m$ is a simply
connected, complete, real analytic K\"a\-hler surface admitting a
negative almost Hermitian structure satisfying the first Gray
condition. Then $M$ is a product of two Riemannian surfaces or
$M=(\Bbb C^2,can)$ with standard flat K\"ahler metric.}
\medskip
{\it Proof.}  We use the classification result of Szabo
(Th.4.5,p.103 in [Sz]). Since $\m$ is a semisymmetric space it is
a direct product of symmetric spaces and Riemannian surfaces,
3-dimensional spacesd which are hyperbolically foliated on an
everywhere dense open subset and k-dimensional spaces which are
parabolically foliated on an open dense subset. From Lemma 3 and
the fact that $dim H^1_3$ , $dim PF^n_{n-2}\ge3$ it follows that
$M$ is symmetric hence space form or is a product of Riemannian
surfaces (note that $P^4_2$ can not be K\"ahler since $dim S=1$
see [Sz]) . The space form of nonzero holomorphic curvature does
not admit a negative almost Hermitian structure satisfying the
first Gray condition (see Lemma 2). Hence $M$ is a product of
Riemannian surfaces or $M=(\Bbb C^2,can)$ with standard flat
K\"ahler metric $can$.
\medskip
{\bf Corollary.}  {\it  Let us assume that $\m$ is a  complete,
real analytic K\"ahler surface admitting a negative almost
Hermitian structure satisfying the first Gray condition. Then $M$
is locally a product of two Riemannian surfaces or $M=\Bbb
C^2\slash \Gamma$.}
\medskip
{\it Remark.}  Note that in the case of product of Riemannian
surfaces the opposite almost Hermitian structure $\bJ$ is
K\"ahler.  In the case $M=(\Bbb C^2,can)$  any opposite almost
Hermitian structure $\bJ$ satisfies the first Gray condition.

\bigskip
{\bf Proposition 6.}  {\it Let $\m$ be a compact K\"ahler surface
admitting an opposite Hermitian structure $\bJ$ satisfying the
first Gray condition. Then $\bJ$ is a K\"ahler structure. If $\m$
is additionaly simply connected then $\m$ is a product of
Riemannian surfaces.}
\medskip
{\it Proof.}  Since $\bJ$ is Hermitian we have
$\kappa=\tau-\frac32(|\th|^2+2\delta\th)$ where $\th$ is the Lee
form of $(M,g,\bJ)$(see [G]). Hence if $\bJ$ satisfies $(G_1)$
then $\tau=\kappa$ and $|\th|^2+2\delta\th=0$. Thus
$$\int_M (|\th|^2+2\delta\th)=0$$ and consequently $\int_M |\th|^2=0$
which gives $\th=0$.  Thus $(M,g,\bJ)$ is K\"ahler and the result
follows from the De Rham theorem.$\k$

\bigskip
Now we give examples of K\"ahler surfaces foliated with two
dimensional Euclidean spaces which are not a product of two
Riemannian surfaces. Hence they admit a negative almost Hermitian
structure satisfying the first Gray condition $G_1$, in fact this
structure is Hermitian. These manifolds are of Calabi type and
hence are ambi-K\"ahler. They are not complete.

Let $(\Sigma,h)$ be a compact Riemannian surface with K\"ahler
form $\om$ such that $\frac1{2\pi}\om$ is an inegral form
corresponding to $1\in H^2(\Sigma,\Bbb Z)=\Bbb Z$.  Let $P_k$ be a
$S^1$ bundle over $\Sigma$ corresponding to an integral class
$k\frac1{2\pi}\om$ where $k\in\Bbb N$.  Let $\th$ be a connection
form on an $S^1$-principal fibre bundle $P_k$ such that
$d\th=k\om$.  Let us consider the manifold  $M=\Bbb R_+\times P_k$
where $\Bbb R_+=\{x\in\Bbb R:x>0\}$ with the metric
$g=dt^2+(\frac{2gg'}k)^2\th^2+g^2p^*h$ where $g=g(t)$ is a
function on $\Bbb R$ and $t$ is the natural coordinate on $\Bbb
R_+$. The metric $g$ is K\"ahler and admits a negative Hermitian
structure $I$. The fundamental vector field $\xi$ of the action of
$S^1$ on $P_k$ is a holomorphic Killing vector field for $M$( see
[J-2]). $M$ is a Calabi type manifold and $a+\frac b2=
-4\frac{g''}g$ (see [J-2],[J-3]). Hence $M$ is semi-symmetric if
$g(t)=t$. $M$ is the fibre bundle over $\Sigma$ with totally
geodesic fibers $\Bbb C^*$ with flat metric. Hence $V^0$ is the
distribution tangent to the fibers $\Bbb C^*$. The metric is
$g=dt^2+\frac{4t^2}{k^2}\th^2+t^2p^*h$.

We shall show now that QCH surfaces with nonvanishing Bochner
tensor $W=W^-$ have the property that the condition $R.W=0$
implies $R.R=0$.
\medskip
{\bf Proposition 6.}  {\it  Let $(M,g,J)$ be a K\"ahler QCH
surface with $R. W^-=0$ and $\kappa\ne0$ on $M$.  Then  $(M,g,J)$
is semi-symmetric.}
\medskip
{\it Proof.}  Since $W^-=\frac{\kappa}2(\frac16\Pi-\Phi+\Psi)$ and
$\kappa\ne0$  it follows that $R.W^-=0$ if and only if $R.K=0$
where $K=\frac16\Pi-\Phi+\Psi$.  Note that $\Phi-\Phi'=\Psi-\Psi'$
and $\Phi+\Phi'=\Pi$.  Hence $\n\Phi=-\n\Phi'$ and
$2\n\Phi=\n\Psi-\n\Psi'$.  Consequently  $2\n K= \n\Psi+\n\Psi'$.
Let us write $\om_1=h_J,\om_2=m_J$.  Then $\om=\om_1+\om_2$  and
$\Psi=-\om_1\otimes\om_1, \Psi'=-\om_2\otimes\om_2$.  Thus
$$\gather \n_X\Psi+\n_X\Psi'=\n_X\om_1\otimes\om_1+\om_1\otimes\n_X\om_1+\tag 2.13\\
\n_X\om_2\otimes\om_2+\om_2\otimes\n_X\om_2=\n_X\om_1\otimes(\om_1-\om_2)+(\om_1-\om_2)\otimes\n_X\om_1=\\
\frac12\n_X\om'\otimes\om'+\frac12\om'\otimes\n_X\om'=\frac12\n_X(\om'\otimes\om').\endgather$$
Hence $R.K=R.(\om'\otimes\om')$. Note that $R.(\om'\otimes\om')=0$
if and only if $R.\om'=0$. In fact
$R.(\om'\otimes\om')=R.\om'\otimes\om'+\om'\otimes R.\om'$.  Let
$\om'(X,Y)=0$ then $0=R.\om'(U,W)\otimes\om'(X,Y)+\om'(U,W)\otimes
R.\om'(X,Y)$ for all $U,W$.  Taking $U,W$ such that
$\om'(U,W)\ne0$ we get $R.\om'(X,Y)=0$.  If $\om'(X,Y)\ne 0$ then
$0=R.(\om'\otimes\om')(X,Y,X,Y)=2R.\om'(X,Y)\om'(X,Y)$ hence again
$R.\om'(X,Y)=0$. Note that $R.\om'=0$  if and only if $\bJ$
satisfies the first Gray condition and hence $\kappa=\tau$ and
$R.R=0$.   We give also another proof of this fact.  Note that
(see [J-1])
$$ \gather R.K=a\Pi.K+b\Phi.K+c\Psi.K=
-a\Pi.\Phi+a.\Pi.\Psi+\frac16b\Phi.\Pi-\tag 2.14
\\b\Phi.\Phi+b\Phi.\Psi+\frac16c
\Psi.\Pi-c\Psi.\Phi+c\Psi.\Psi=-a\Pi.\Phi+a\Pi.\Psi\\ -b
\Phi.\Phi+b\Phi.\Psi=-a\Pi.\Phi+a\Pi.\Psi-b\frac12\Pi.\Phi+b\frac12\Pi.\Psi.\endgather$$

Since the tensors $\Pi.\Phi,\Pi.\Psi$ are linearly independent it
follows that $a+\frac12b=0$. $\k$

Now we consider semi-symmetric K\"ahler surfaces foliated by
2-dimensional euclidean space.  Let  $\De=V^0=\{X: R(U,V)X=0
\text{    for all  } U,V\in TM\}$.  Then $\De$ is totally geodesic
foliation.  Let $I$ be defined  $IX=JX$  if $X\in\De$ and $IX=-JX$
if $X\in \Cal E=\De^{\perp}$.  Note that  $R=c\Psi$ with respect
to $\Cal E$ and $\tau=2c$ where $\tau$ is the scalar curvature of
$(M,g,J)$. We have
$$\gather
\n_XR(Y,Z,W,T)=-Xc\om(Y,Z)\om(W,T)-c\n_X\om(Y,Z)\om(W,T)\\-c\om(Y,Z)\n_X\om(W,T)
\endgather $$

where $\om=m_J\in\bigwedge^2\E$ and hence from the Bianchi
identity we obtain
$$-dc\w\om\otimes \om(W,T)-cd\om\otimes\om(W,T)-c\om\w\n_.\om(W,T)=0\tag 2.15$$
Since $\rho=c\om$ where $\rho$ is the Ricci form of $(M,g,J)$ we
get $d\om=-d\ln c\om$.  Hence we get from (2.15)
$$\om\w\n_.\om(W,T)=0.$$  Consequently  $\n_X\om=0$ for $X\in\De$.
Note that the K\"ahler form corresponding to $J$ is
$\Om=\om_1+\om$ and corresponding to $I$ is $\Om_1=\om_1-\om$ and
consequently $\n_XI=0$ for $X\in\De$ where
$\om_1=h_J\in\bigwedge^2\De$. We also have $d\Om_1=2d\om_1=2d\ln
c\w\om=2\th\w\Om_1$ where $th(X)=g(\n\ln c_{|\De},X)$ is the Lee
form of $I$.  Let us assume that $\zeta\in\G(\De)$.  Then
$$L_{\zeta}\rho=\zeta\lrcorner d\rho+d(\zeta\lrcorner\rho)=0.$$
We also have $\n_{\zeta}\rho=\zeta c\om$ since $\n_{\zeta}\om=0$.
Thus
$(L_{\zeta}-\n_{\zeta}).\rho(X,Y)=-\n\zeta.\rho(X,Y)=c\om(\n_X\zeta,Y)+c\om(X,\n_Y\zeta).$
Hence $\om(\n_X\zeta,Y)+\om(X,\n_Y\zeta)=-\th(\zeta)\om(X,Y)$. Let
us assume that $X,Y\in\G(\Cal E)$. Then
$$g(\n_XJ\zeta,Y)-g(X,\n_YJ\zeta)=-\th(\zeta)\om(X,Y)$$
and for any $\xi\in\G(\De )$ we obtain
$$g(\n_X\xi,Y)-g(X,\n_Y\xi)=\th(J\xi)\om(X,Y).\tag2.16$$
Note that since $\De$ is totally geodesic the foliation $\De$ is
holomorphic if and only if  $g(\_{JX}\xi,Y)=g(J\n_X\xi,Y)$ for any
$X,Y\in\Cal E$. Hence if $\De$ is holomorphic we get
$$\gather g(\n_{JX}\xi,Y)-g(JX,\n_Y\xi)=-\th(J\xi)g(X,Y),\\
g(\n_{X}J\xi,Y)+g(X,\n_YJ\xi)=-\th(J\xi)g(X,Y)\\
g(\n_{X}\zeta,Y)+g(X,\n_Y\zeta)=-\th(\zeta)g(X,Y)\tag
2.17\endgather$$ for any $\zeta\in \De$.  From (2.16) and (2.17)
we obtain
$$2g(\n_X\zeta,Y)=-J\th(\zeta)\om(X,Y)-\th(\zeta)g(X,Y).\tag 2.18$$
Hence
$$2(\n_XY)_{|\De}=\th^\sharp g(X,Y)+J\th^\sharp \om(X,Y)\tag 2.18$$ for all
$X,Y\in\G(\E)$. From (2.17) it follows that if $\De$ is
holomorphic then it is conformal.  On the other hand let us assume
that $\De$ is conformal i.e.
$g(\n_{X}\zeta,Y)+g(X,\n_Y\zeta)=\phi(\zeta)g(X,Y)$ for some form
$\phi\in\bigwedge^1\De$.  Then
$$2g(\n_X\zeta,Y)=-J\th(\zeta)\om(X,Y)+\phi(\zeta)g(X,Y).\tag 2.19$$
From (2.19) it is clear that $g(\n_{JX}\zeta,JY)=g(\n_X\zeta,Y)$
for all $X,Y\in\E$ and thus $\n_{JX}\zeta-J\n_X\zeta\in\De$ which
means that $\De$ is holomorphic.  It follows that $\phi=-\th$.

Thus we get

{\bf Theorem} {\it  Let $(M,g,J)$ be a K\"ahler semi-symmetric
surface foliated by two dimensional Euclidean space.  Then the
following conditions are equivalent:

(a) $\De$ is a holomorphic foliation

(b) $\De$ is a conformal foliation

(c)  $2g(\n_X\zeta,Y)=-J\th(\zeta)\om(X,Y)-\th(\zeta)g(X,Y)$ for
all $\zeta\in\G(\De)$ and all $X,Y\in\G(\E)$

(d)  The almost hermitian structure $I$ is Hermitian.}

\medskip
{\it Proof.}  The equivalence of conditions (a),(b),(c) we have
proved above.  We shall show that conditions (a) and (d) are
equivalent.  Note that if $Y\in\G(\De)$ then $IY=JY$ and
consequently $\n_XIY=2J(\n_XY_{|\E}$.  Similarly if $Y\in\G(\E)$
then $\n_XIY=-2J(\n_XY_{|\De}$.  Let us assume that (a) holds.
Then (c) holds.  To prove that $I$ is integrable we have to show
that $\n_XIY=\n_{IX}I IY$ for $X\in\G(\E)$ and $Y\in\frak X(M)$.
Let $\{E_1,E_2,E_3,E_4\}$ be a local orthonormal basis of $(M,g)$
such that $E_1,E_2$ span $\De$ and $E_3,E_4$ span $\E$. We assume
further that $JE_1=IE_1=E_2, JE_3=-IE_3=E_4$. We have
$\n_{E_3}IE_1=2J(\n_{E_3}E_1)_{|\E}$ and
$\n_{IE_3}IIE_4=-\n_{E_4}IE_2=-2J(\n_{E_4}E_2)_{|\E}$.  Using (c)
it is cler that $\n_{E_3}IE_1=\n_{IE_3}IIE_4$.  If $X,Y\in\G(E)$
then  $\n_XIY=-2J(\n_XY)_{|\De}$ and
$\n_{IX}IIY=-2J(\n_{JX}JY)_{|\De}$.  From (c) it is clear that
$\n_XIY=\n_{IX}IIY$ also in this case.   On the other hand if $I$
is integrable then $\n_XY_{|\De}=\n_{JX}JY)_{|\De}$ for all
$X,Y\in\G(\E)$ which is equivalent to $\De$ being holomorphic.
$\k$

Since the Weyl tensor $W^-$ is degenerate the form $d\th$ is
self-dual and $d\th(J,J)=-d\th$ if $I$ is degenerate.

{\bf Proposition. } {\it  Let $(M,g,J)$ be a K\"ahler
semi-symmetric surface foliated by two dimensional Euclidean
space.  Then the following conditions are equivalent:
 (a)  $d\th=0$

 (b) $J\th^\sharp$ is a Hamiltonian vector field for $(M,g,J)$

 (c) $Y|\th|^2=0$ for $Y\in\E$

 Moreover if $I$ is integrable every of these condition is
 equivalent to

 (e) $I$ is locally conformally K\"ahler}

\medskip

{\it Proof.}  Note that $\om_1=|\th|^{-2}\th\w J\th$ and
$d\om_1=\th\w\om$.  Hence $L_{J\th}\om_1=-d\th$. It is clear that
$L_{J\th}\om=0$.  Hence $L_{J\th}\Om=-d\th=L_{J\th}\Om_1$. Hence
(a) is equivalent to (b).  It is easy to see that $d\th(X,Y)=0$ if
$X,Y$ are both in $\De$ or $\E$.  Now
$d\th(Y,\xi)=\frac12Y|\th|^2$ for $Y\in \E$.  Since
$d\th(J,J)=-d\th$ this proves that (a) is equivalent to (c).$\k$

In the following theorem we use description of Calabi type K\"aler
surfaces from the paper [A-C-G].

\medskip
 {\bf Theorem. } {\it Let $(M,g,J)$ be a
K\"ahler surface admitting opposite Hermitian  structure $I$
satisfying the first Gray condition $G_1$ which is locally
conformally K\"ahler. Then locally

$$ g=zg_{\Sigma}+\frac1{Cz}dz^2+Cz(dt+\alpha)^2\tag 2.19$$
where $(\Sigma,g_{\Sigma})$ is a Riemannian surface with area form
$\om_{\Sigma}$ and $d\alpha=\om_{\Sigma}$ or $(M,g,J)$ is a
product of Riemannian surfaces or a space form with zero
holomorphic sectional curvature. The K\"ahler form of $(M,g,J)$ is
$\Om=z\om_{\Sigma}+dz\w(dt+\alpha)$. }

\medskip
{\it Proof.} First assume that $(M,g,J)$ is a K\"ahler
semi-symmetric surface foliated by 2-dimensional Euclidean space.
 Hence it follows
that $\De$ is totally geodesic homothetic foliation. Such
foliations were classified locally  in [Ch-N]. Thus $(M,g,J)$ is a
K\"ahler surface of Calabi type.  From [A-C-G] it follows that
$\kappa=\tau$ if $V(z)=Cz^2$ where $$ g=zg_{\Sigma}+\frac
z{V(z)}dz^2+\frac {V(z)}z (dt+\alpha)^2$$ is a general Calabi type
metric which is not a K\"ahler product. For the general case let
us note that QCH K\"ahler surface for which the structure $I$ is
Hermitian and locally conformally K\"ahler are of Calabi type or
are orthotoric surfaces or $W=0$ (see [J-4]). Semi-symmetric
surfaces with $W=0$ are products of Riemannian surfuces of
constant opposite scalar curvatures (see [B]). One can easily
check that orthotoric surface can be semi-symmetric only if $W=0$
which finishes the proof.   $\k$

Note that the examples with the metric
$g=dt^2+(\frac{2gg'}k)^2\th^2+g^2p^*h$ given before are the
special kind of (2.19) and these classes of manifolds coincide
locally.

The paper was supported by NCN grant 2011/01/B/ST1/02643.

\bigskip
\centerline{\bf References.}

\par
\medskip
\cite{A-C-G} V. Apostolov,D.M.J. Calderbank, P. Gauduchon {\it The
geometry of weakly self-dual K\"ahler surfaces}  Compos. Math.
135, 279-322, (2003)
\par
\medskip
\cite{B}  Bryant R.  {\it Bochner-K\"ahler metrics} J. Amer. Math.
Soc.14 (2001) , 623-715.
\par
\medskip
[Ch-N] S.G.Chiossi and P-A. Nagy {\it Complex homothetic
foliations on K\"ahler manifolds}  Bull. London Math. Soc. 44
(2012) 113-124.
\par
\medskip
\cite{D} A. Derdzi\'nski, {\it Self-dual K\"ahler manifolds and
Einstein manifolds of dimension four }, Compos. Math.
49,(1983),405-433

\par
\medskip
\cite{G-M-1} G.Ganchev, V. Mihova {\it K\"ahler manifolds of
quasi-constant holomorphic sectional curvatures}, Cent. Eur. J.
Math. 6(1),(2008), 43-75.
\par
\medskip
\cite{G-M-2} G.Ganchev, V. Mihova {\it Warped product K\"ahler
manifolds and Bochner-K\"ahler metrics}, J. Geom. Phys. 58(2008),
803-824.
\par
\medskip
\cite{G} P. Gauduchon {\it La 1-forme de torsion d'une variete
hermitienne compacte} Math. Ann. 267, (1984),495-518.
\par
\medskip
\cite{J-1} W. Jelonek,{\it Compact holomorphically pseudosymmetric
K\"ahler manifolds} Coll. Math.117,(2009),No.2,243-249.
\par
\medskip
\cite {J-2} W.Jelonek {\it K\"ahler manifolds with quasi-constant
holomorphic curvature}, Ann. Glob. Anal. and Geom, vol.36, p.
143-159,( 2009)
\par
\medskip
\cite{J-3} W. Jelonek, {\it Holomorphically pseudosymmetric
K\"ahler metrics on $\Bbb{CP}^n$} Coll.
Math.127,(2012),No.1,127-131.
\par
\medskip
\cite{J-4} W. Jelonek {\it K\"ahler surfaces with quasi-constant
holomorphic curvature} arxiv
\par
\medskip
\cite{L} U. Lumiste {\it Semisymmetric curvature operators and
Riemannian 4-spaces elementary classified}  Algebras, Groups and
Geometries 13, 371-388 (1996)
\par
\medskip
\cite{O} Z. Olszak, {\it Bochner flat K\"ahlerian manifolds with a
certain condition on the Ricci tensor} Simon Stevin 63,
(1989),295-303
\par
\medskip
\cite{K-N} S. Kobayashi and K. Nomizu {\it Foundations of
Differential Geometry}, vol.2, Interscience, New York  1963
\par
\medskip
\cite{Sz} Z. I. Szabo {\it Structure theorems on Riemannian spaces
satisfying $R(X,Y).R$ $=0$ II.}  Global version  Geometriae
Dedicata {\bf 19} (1985), 65-108.

\par
\medskip
Institute of Mathematics

Cracow University of Technology

Warszawska 24

31-155      Krak\'ow,  POLAND.

 E-mail address: wjelon\@pk.edu.pl
\bigskip

\enddocument